\documentclass[11pt]{article}
\usepackage{amstext,amssymb,amsmath,amsbsy}

\usepackage{hyperref}
\usepackage{amscd}
\usepackage{amsfonts}
\usepackage{indentfirst}
\usepackage{verbatim}
\usepackage{amsmath}
\usepackage{amsthm}
\usepackage{enumerate}
\usepackage{graphicx}
\usepackage{subfig}
\usepackage{color}
\usepackage[OT1]{fontenc}
\usepackage[latin1]{inputenc}
\usepackage[english]{babel}
\usepackage{amssymb}
\newtheorem{theorem}{Theorem}

\newtheorem{remark}{Remark}

\numberwithin{equation}{section}

\newcommand{\proofend}{\hfill $\Box$ }

\newcommand{\dsp}{\displaystyle}

\newcommand{\supp}{\operatorname{supp}}

\newcommand{\dive}{\operatorname{div}}

\newcommand{\eps}{\varepsilon}

\newcommand{\epss}{s}

\newcommand{\loc}{_{loc}}

\newcommand{\mC}{\mathbb{C}}

\newcommand{\mR}{\mathbb{R}}

\newcommand{\AC}{AC}
\newcommand{\AD}{AD}
\newcommand{\BC}{BC}
\newcommand{\BD}{BD}
\newcommand{\EC}{EC}
\newcommand{\ED}{ED}
\newcommand{\FC}{FC}
\newcommand{\FD}{FD}

\title{Superlensing using complementary media}

\author{Hoai-Minh Nguyen \footnote{EPFL SB MATHAA CAMA, Station 8,  CH-1015 Lausanne, hoai-minh.nguyen@epfl.ch} \footnote{School of Mathematics, University of Minnesota, MN, 55455, hmnguyen@math.umn.edu} \footnote{The research is supported by NSF grant DMS-1201370 and by the Alfred P. Sloan Foundation.}}

\begin{document}

\maketitle 

\begin{abstract}
This paper studies magnifying superlens using complementary media. Superlensing using complementary media was suggested by Veselago in \cite{Veselago} and innovated by Nicorovici et al. in \cite{NicoroviciMcPhedranMilton94} and Pendry in \cite{PendryNegative}.  The study of this problem is difficult due to two facts. Firstly, this problem is unstable since the equations describing the phenomena  have sign changing coefficients; hence the ellipticity is lost. Secondly, the phenomena associated are localized resonant, i.e., the field explodes in some regions and remains bounded in some others. This makes the problem difficult to analyse.  In this paper, we develop the technique of removing of localized singularity introduced in \cite{Ng-Negative-Cloaking} and make use of the reflecting technique  in \cite{Ng-Complementary} to overcome these two difficulties.  More precisely, we suggest a class of  lenses which has root from \cite{NicoroviciMcPhedranMilton94} and \cite{PendryRamakrishna} and inspired from \cite{Ng-Negative-Cloaking} and give a proof of superlensing for this class.  To our knowledge, this is the first rigorous proof  on the  magnification of an  {\bf arbitrary inhomogeneous} object using complementary media. 
\end{abstract}

\section{Introduction}

Negative index materials (NIMs) were first investigated theoretically by Veselago in \cite{Veselago} and innovated by Nicorovici et al. in \cite{NicoroviciMcPhedranMilton94} and Pendry in \cite{PendryNegative}. The  existence of such materials was confirmed by Shelby et al. in \cite{ShelbySmithSchultz}. NIMs have been intensively  studied recently thanks to its many applications and surprising properties. One of the appealing ones is superlensing. The construction of a slab superlens using NIMs was suggested by Veselago in \cite{Veselago} via the ray theory. Later, this was developed by  Nicorovici et al.  in \cite{NicoroviciMcPhedranMilton94}  and Pendry in \cite{PendryNegative}. In \cite{NicoroviciMcPhedranMilton94} the authors studied a cylindrical lens in the two dimensional  quasistatic regime, and  in  \cite{PendryNegative} the author studied the Veselago slab in the finite frequency one. These works have been developed further, see, e.g.,  in  \cite{MiltonNicoroviciMcPhedranPodolskiy, PendryCylindricalLenses, PendryRamakrishna-1, PendryRamakrishna0, PendryRamakrishna} where  cylindrical and spherical superlenses were investigated. The reader can find an interesting review and many recent results on  superlensing using complementary media in \cite{MiltonNicoroviciMcPhedranPodolskiy}. 

\medskip
The study of superlensing has been  concentrating a lot on  the image of dipoles in homogeneous media see \cite{MiltonNicoroviciMcPhedranPodolskiy, PendryNegative,PendryCylindricalLenses, PendryRamakrishna-1, PendryRamakrishna0, PendryRamakrishna}. There are a few  works devoted to the image of an object. It seems for us that \cite{NicoroviciMcPhedranMilton94}, in which the authors gave a proof on the magnification of a constant material disk,  is the only work  in this direction. Even though, the methods in the papers mentioned above can be used to obtain the magnification of radial objects having constant materials in two or three dimensions, the magnification of an arbitrary inhomogeneous object is out of scope of these methods, which are strongly based on the separation of variables. Let us mention two difficulties related to the study of this problem. Firstly, this problem is unstable since the equations describing the phenomena  have sign changing coefficients; hence the ellipticity is lost. Secondly, the phenomena associated are localized resonant, i.e., the field explodes in some regions and remains bounded in some others. This makes the problem difficult to analyze. 

\medskip
In this paper, we study magnifying superlens using complementary media.  More precisely, given $m>1$ the magnification, 
we suggest a class of lenses, which has root from \cite{NicoroviciMcPhedranMilton94} and \cite{PendryRamakrishna} and inspired from \cite{Ng-Negative-Cloaking}, and  show that one can magnify $m$ times an {\bf arbitrary inhomogeneous} object in the quasistatic and finite frequency regimes using a lens in this class.   To overcome these difficulties mentioned above,  we develop the technique of removing localized singularity introduced in \cite{Ng-Negative-Cloaking}, and make use of the reflecting technique in \cite{Ng-Complementary}. To our knowledge, these results of this paper are new even in the two dimensional quasistatic regime. 

\medskip
Let us describe how to magnify the region $B_{r_0}$ for some $r_0> 0$ in which the medium is characterized by a matrix-valued function $a$ and a real function $\sigma$ using complementary media. Here and in what follows given $r > 0$, $B_r$ denotes the ball centered at the origin of radius $r$ in $\mR^d$ ($d = 2$ or 3). The assumption on the geometry of the object by all means imposes no restriction since any region can be placed in such  a ball provided that the radius and the origin are appropriately chosen. We first concentrate on  the quasistatic regime.  The idea suggested in \cite{NicoroviciMcPhedranMilton94, PendryCylindricalLenses, PendryRamakrishna} is to put a lens in $B_{r_2} \setminus B_{r_0}$ whose medium is  characterized by matrix $-b$  with $r_2^2/ r_0^2 = m$. Here  $b = I$, the identity matrix,  in two dimensions and $b = \big(r_2^2/|x|^2 \big)I$ in three dimensions. 

\medskip In this paper,  we slightly change  the strategy discussed above and take into account the suggestion in \cite{Ng-Negative-Cloaking}.  Our lens contains two parts. The first one is given by 
\begin{equation}\label{first-part}
- (r_2^2/ |x|^2)^{d-2} I \quad \mbox{ in } B_{r_2} \setminus B_{r_1}
\end{equation} 
and the second one is the matrix 
\begin{equation}\label{second-part}
m^{d-2} I \quad \mbox{ in } B_{r_1} \setminus B_{r_0}. 
\end{equation}
Here 
\begin{equation}\label{def-r1}
r_1 = m^{1/4} r_0
\end{equation}
and 
\begin{equation}\label{def-r2}
r_2 =   \sqrt{m} r_1 =  m^{3/4} r_0. 
\end{equation}
Set
\begin{equation}\label{def-r3}
r_3 := r_2^2/ r_1 = m^{5/4} r_0. 
\end{equation}
It is clear that 
\begin{equation}\label{pro-m}
m = r_2^2/ r_1^2. 
\end{equation}
We will give some comments on this construction later. 

\medskip
Since materials have some loss, the correct approach is to allow some loss in the medium and investigate the limit as the loss goes to 0. With the loss, the medium is characterized by $s_\delta A$, 
where 
\begin{equation}\label{defA}
A = \left\{ \begin{array}{cl}  \big( r_2^2 / |x|^2 \big)^{d-2} I & \mbox{ in } B_{r_2} \setminus B_{r_1}, \\[6pt]
m^{d-2} I & \mbox{ in } B_{r_1} \setminus B_{r_0}, \\[6pt]
a & \mbox{ in } B_{r_0},\\[6pt]
I  & \mbox{ otherwise}, 
\end{array} \right. 
\end{equation} and 
\begin{equation}\label{def-ss}
\epss_\delta = \left\{ \begin{array}{cl} -1 + i \delta  & \mbox{ in } B_{r_2} \setminus B_{r_1}, \\[6pt]
1 & \mbox{ otherwise}. 
\end{array} \right. 
\end{equation}
Physically, the imaginary part of $\epss_\delta A$ is the loss of the medium (more precisely the loss of  the medium in 
$B_{r_2} \setminus B_{r_1}$).
In what follows, we assume that
\begin{equation}\label{cond-a}
\frac{1}{\Lambda} |\xi|^{2} \le \langle a (x) \xi, \xi \rangle \le \Lambda |\xi|^{2} \quad \forall \, \xi \in \mR^{d}, \mbox{ for a.e. } x \in B_{R_{3}},    
\end{equation}
for some constant $\Lambda \ge 1$. 

\medskip
We next make some comments on the construction. We first note that $-(r_2^2/|x|^2)^{d-2} I$ in $B_{r_2} \setminus B_{r_1}$ and $I$ in  $B_{r_3} \setminus B_{r_2}$ are complementary or more precisely reflecting complementary via the Kelvin transform $F: B_{r_2} \to \mR^d \setminus \bar B_{r_2}$ w.r.t. $\partial B_{r_2}$, i.e, 
\begin{equation}\label{F-b}
F(x) = r_2^2 x/ |x|^2 \quad \mbox{ and } \quad   F_*A = I \mbox{ in } B_{r_3} \setminus B_{r_2}.
\end{equation}
(see \cite{Ng-Complementary} for the definition of reflecting complementary media and their properties). 
Here
\begin{equation}\label{change1}
T_*M(y) = \frac{D T  (x)  M(x) D T ^{T}(x)}{J(x)} \mbox{ where }  x =T^{-1}(y) \mbox{ and } J(x) = |\det D T(x)|, 
\end{equation}
for a diffeomorphism $T$ and a matrix $M$. Given $r_1$,  
the choice of $r_2$ follows from \eqref{pro-m} since a superlens of  $m$ times magnification is considered as in \cite{NicoroviciMcPhedranMilton94, PendryCylindricalLenses, PendryRamakrishna} (see also \eqref{def-hatA} and Theorem~\ref{thm1}). The choice of $r_1$ and $A$ in $B_{r_1} \setminus B_{r_0}$ are inspired from \cite{Ng-Complementary, Ng-Negative-Cloaking} as follows. Let $G: \mR^d \setminus \bar B_{r_3} \to B_{r_3} \setminus \{0\}$  be the Kelvin transform w.r.t. $\partial B_{r_3}$, i.e., 
\begin{equation}\label{G-b}
G(x) = r_3^2 x / |x|^2. 
\end{equation}
Then $G \circ F: B_{r_1} \to B_{r_3}$ satisfies
\begin{equation}\label{GF}
G\circ F (x) = m x  \mbox{ in } B_{r_1}.  
\end{equation} 
This implies, since $A= m^{d-2} I$ in $B_{r_1} \setminus B_{r_0}$.  
\begin{equation}\label{the-choice}
G_*F_* A  = I \mbox{ in } B_{r_3} \setminus B_{r_{*}}.  
\end{equation}
Here 
\begin{equation}\label{def-r*}
r_* := m r_0 = \sqrt{r_2 r_3}= \sqrt{r_2^3/r_1}. 
\end{equation} 
In the last identity, we use the fact that  $r_3 = r_2^2/ r_1$. 
Using \eqref{pro-m} and \eqref{def-r*}, we derive the formula for $r_1$ and $r_2$ as in \eqref{def-r1} and \eqref{def-r2}.  
The choice of $A$ in $B_{r_1} \setminus B_{r_0}$  follows from \eqref{the-choice}. 

\medskip
In the finite frequency regime, the medium is also characterized by $s_\delta\Sigma$ where 
\begin{equation}\label{id1}
\Sigma = \left\{\begin{array}{cl} \dsp \big(r_2^2/ |x|^{2} \big)^d & \mbox{ if } x \in B_{r_{2}} \setminus B_{r_1}, \\[6pt]
\dsp m^{d} & \mbox{ if } x \in B_{r_{1}} \setminus B_{r_0}, \\[6pt]
\sigma & \mbox{ in } B_{r_0}, \\[6pt]
1 & \mbox{ otherwise}. 
\end{array}\right.
\end{equation} 
The construction of $\Sigma$ for the lens is given in $B_{r_2} \setminus B_{r_0}$. This construction is based on the requirements
\begin{equation}\label{choice-sigma}
F_* \Sigma = 1 \mbox{ in } B_{r_3} \setminus B_{r_2} \quad \mbox{ and } \quad G_*F_*\Sigma = 1 \mbox{ in } B_{r_3} \setminus B_{r_*}.  
\end{equation}
Here
\begin{equation}\label{change2}
T_*h(y) = \frac{h(x) }{J(x)} \mbox{ where }  x =T^{-1}(y) \mbox{ and } J(x) = |\det D T(x)|, 
\end{equation}
for a diffeomorphism $T$ and a function $h$.
These requirements  are not easy to predict but follow naturally from the study of reflecting complementary media in  \cite{Ng-Complementary}.  
We will assume that 
\begin{equation}\label{assumption-sigma}
1/ \Lambda \le \sigma(x) \le \Lambda \quad \mbox{ for a.e. } x \in B_{r_0},   
\end{equation}
for some $\Lambda \ge 1$. 


\medskip

This paper deals with the bounded setting equipped  the zero Dirichlet boundary condition. Let $k \ge 0$ and $\Omega$ be a smooth open subset of $\mR^d$ ($d=2, \, 3$) such that $B_{r_3} \subset \Omega$. Given $f \in L^{2}(\Omega)$, let  $u_\delta, \, u \in H^1_{0}(\Omega)$ be resp. the unique solution (the well-posedness follows from \eqref{WP} and \eqref{WP1} below) to 
\begin{equation}\label{eq-uu-delta}
\dive (\epss_\delta A \nabla u_\delta) +  s_\delta k^2 \Sigma = f \mbox{ in } \Omega, 
\end{equation}
and 
\begin{equation}\label{eq-uu}
\dive(\hat A \nabla u) + k^2 \hat \Sigma u = f \mbox{ in } \Omega. 
\end{equation}
Here 
\begin{equation}\label{def-hatA}
\hat A,\hat \Sigma = \left\{ \begin{array}{cl} 
m^{2-d} a (x/m), m^{-d} \sigma (x/m)  & \mbox{ in } B_{m r_0},\\[6pt]
I,1  & \mbox{ otherwise}.
\end{array} \right.
\end{equation}

When $k>0$, we will assume in addition that, as in \cite{Ng-Complementary},  
\begin{equation}\label{WP}
 \eqref{eq-uu} \mbox{ is well-posed in } H^1_0(\Omega)
\end{equation}
and  
\begin{equation}\label{WP1}
\mbox{the equation } \Delta v + k^2 v = 0 \mbox{ in } \Omega \setminus B_{r_2} \mbox{ has only zero solution in } H^1_0(\Omega \setminus B_{r_2}). 
\end{equation}

\medskip
Here is one of the two main results of this paper (the second one is Theorem~\ref{thm2} in Section~\ref{sect-comments}). 
\begin{theorem}\label{thm1} Let  $d=2, \, 3$, $f \in L^{2}(\Omega)$ with $\supp f \subset \Omega \setminus B_{r_{3}}$ and let  $u, \; u_\delta  \in  H^1_{0}(\Omega)$ be the unique solutions to \eqref{eq-uu-delta} and \eqref{eq-uu} resp. We have 
\begin{equation}\label{key-point}
u_{\delta} \to u \mbox{ weakly in } H^{1} (\Omega \setminus B_{r_3}) \mbox{ as } \delta \to 0. 
\end{equation}
\end{theorem}

For an observer outside $B_{r_3}$, the object $(a, \sigma)$ in $B_{r_0}$ would act like 
\begin{equation*}
\big( m^{2-d} a (x/m), m^{-d} \sigma(x/m) \big) \mbox{ in } B_{m r_0}
\end{equation*}
by \eqref{key-point}: one has a superlens whose the magnification is $m$.

\medskip  The key ingredient of the proof of Theorem~\ref{thm1} is the removing of localized singularity technique which is introduced in \cite{Ng-Negative-Cloaking} to study cloaking using complementary media. The reflecting technique, which is presented in \cite{Ng-Complementary} also plays an important role in our analysis. In \cite{Ng-CALR}, these techniques will be developed for the context of cloaking due to anomalous localized resonance. To make use of these techniques, we require that $A = m^{d-2} I$ and $\Sigma = m^{d-2}$ in $
B_{r_1} \setminus B_{r_0}$ (which is the second part of our lens construction). Indeed, in the proof we use interpolation inequalities in which the conditions $r_* \le \sqrt{r_2 r_3}$, $ G_*F_A = I$, $G_*F_*\Sigma = 1$ are required, see,  e.g., \eqref{re-u2} and \eqref{ineq-Holder}.  It was argued in \cite{MiltonNicoroviciMcPhedranPodolskiy} that in the two dimensional quasistatic regime,  to be successfully imaged, a conducting object has to be placed in the circle $B_r$ with $r \le \sqrt{r_1^3/r_2}$. In our notations, it is required that $r_1 \ge m^{1/4} r_0$; hence the layer $B_{r_1} \setminus B_{r_0}$ might be necessary. Nevertheless, we do not know how to prove or disprove the necessity of this layer.  

\medskip
It was showed in \cite[Theorem 1]{Ng-Complementary} that 
\eqref{key-point} holds if $\|u_{\delta}\|_{H^{1}}$ remains bounded (this is equivalent to the compatibility condition on $f$ in \cite[Definition 2]{Ng-Complementary}). The goal of this paper is to show that \eqref{key-point} holds without the compatibility assumption. It is clear that the localized resonance appears if the compatibility does not hold. 
The localized resonance appears in this situation would be anomalous one whose concept is introduced in \cite{MiltonNicorovici} since  it seems that the boundary of the resonant regions would vary with the position of the source, and their boundary do not coincide with any discontinuity in moduli. We do not verify this property in this work.  We note that there are plasmonic structures for which either localized resonance or else completely resonance takes places whenever resonance appears see \cite{MinhLoc}. The localized resonance is related to the geometry of the problem. 

\medskip
The lens  in the region $B_{r_2} \setminus B_{r_1}$ discussed above is given by $I$ in two dimensions and $(r_2^2/ |x|^2) I$ in three dimensions. The construction in three dimensions from \cite{PendryRamakrishna0, PendryRamakrishna} is more involved than the one in two dimensions and based on the search of  isotropic radial forms. In Section~\ref{sect-comments}, we will extend this construction to a class of lenses containing anisotropic ones (Theorem~\ref{thm2}). In particular, we will point out a construction for which $r_3$ can be arbitrary close to $m r_0$ (see Remark~\ref{thickness}).  
This extension is based on the study of reflecting complementary media in \cite{Ng-Complementary}.
The concept of complementary media was originally suggested in \cite{PendryRamakrishna-1, PendryRamakrishna0} (see also \cite{LaiChenZhangChanComplementary, NicoroviciMcPhedranMilton94, PendryNegative, PendryRamakrishna}), where various examples were mentioned , and played an important role in the study of NIMs.  In \cite{Ng-Complementary}, the author provides a precise definition of a class of complementary media, reflecting complementary media, generated by reflections and investigates the properties of this class.

\medskip 
The paper is organized as follows. The proof of Theorem~\ref{thm1} will be given in Section~\ref{sect-thm1}. Theorem~\ref{thm2}, a generalization of Theorem~\ref{thm1} which allows anisotropic lenses, will be given in Section~\ref{sect-comments}.

\section{Proof of Theorem~\ref{thm1}}\label{sect-thm1} 

This section is devoted to the proof of Theorem~\ref{thm1}. We first present the proof in the two dimensional quasistatic case (Section~\ref{part1}). We will profit the notational ease in this case to present clearly the ideas of the proof. The proof in the three dimensional quasistatic case is  briefly sketched in Section~\ref{part2}. In Section~\ref{part3}, we consider the finite frequency case. The proof in this case is similar to the one in the quasistatic one though more involved, in particular, for low modes.  

\subsection{Proof of Theorem~\ref{thm1} in the two dimensional quasistatic regime}\label{part1}

In this section, $k=0$ and $d=2$. Multiplying \eqref{eq-uu-delta} by $\bar u_\delta$ (the conjugate of $u_\delta$), integrating in $\Omega$, and using the fact that $u_\delta = 0$ on $\partial \Omega$, we have 
\begin{equation*}
\int_{\Omega} s_\delta \langle A \nabla u_\delta, \nabla u_\delta \rangle = - \int_{\Omega} f \bar u_\delta. 
\end{equation*}
Considering first the imaginary part and then the real part, we obtain, by \eqref{cond-a}, 
\begin{equation}\label{est-F}
\| u_\delta\|_{H^1(\Omega)}^2 \le \frac{C}{\delta} \| u_\delta\|_{L^2(\Omega \setminus B_{r_3})} \| f\|_{L^2}.  
\end{equation}
Here and in what follows in the proof, $C$ denotes a positive constant independent of $\delta$ and $f$. 

\medskip
As in \cite{Ng-Complementary, Ng-Negative-Cloaking}, let $u_{1, \delta} \in H^1_{\loc}(\mR^d \setminus B_{r_2})$ be the reflection of $u_\delta$  through $\partial B_{r_2}$ by $F$, i.e., 
\begin{equation*}
u_{1, \delta} = u_\delta \circ F^{-1},  
\end{equation*}
and let $u_{2, \delta} \in H^1(B_{r_3})$ be the reflection of $u_{1, \delta}$  through $\partial B_{r_3}$ by $G$, i.e., 
\begin{equation*}
u_{2, \delta} = u_{1, \delta} \circ G^{-1} = u_\delta \circ F^{-1} \circ G^{-1}.  
\end{equation*}
We recall that $F$ and $G$ are given in \eqref{F-b} and \eqref{G-b}. Since $G\circ F (x) = \big(r_3^2/ r_2^2 \big) x$, it follows from \eqref{change1} that 
\begin{equation*}
\hat A = G_*F_*A \mbox{ in } B_{r_3}. 
\end{equation*}
Applying \cite[Lemma 2]{Ng-Complementary}, we have 
\begin{equation}\label{Delta-u1}
\Delta u_{1, \delta} = 0 \quad  \mbox{ in } B_{r_3} \setminus B_{r_2} 
\end{equation}
and 
\begin{equation}\label{Delta-uu2}
\dive ( \hat A \nabla u_{2, \delta} ) = 0 \quad  \mbox{ in } B_{r_3}.  
\end{equation}
We derive from \eqref{Delta-uu2} that 
\begin{equation}\label{Delta-u2}
\Delta u_{2, \delta} = 0 \mbox{ in } B_{r_3} \setminus B_{r_*}. 
\end{equation}

From the transmission conditions on $\partial B_{r_2}$, we have
\begin{equation}\label{transmission-r2}
u_{1, \delta} = u_\delta \quad \mbox{ and } \quad  (1 - i \delta)\partial_r u_{1, \delta} = \partial_r u_{\delta} \Big|_{\rm ext} \quad \mbox{ on } \partial B_{r_2}
\end{equation}
and, from the transmission conditions on $\partial B_{r_3}$, we obtain
\begin{equation}\label{transmission-r3}
u_{2, \delta} = u_{1, \delta} \quad \mbox{ and } \quad  \partial_{r} u_{2, \delta}  = (1 - i \delta)\partial_r u_{1, \delta} \Big|_{\rm int} \quad \mbox{ on } \partial B_{r_3}. 
\end{equation}
Since $\Delta u_\delta = \Delta u_{1, \delta} = 0$ in $B_{r_3} \setminus B_{r_2}$, by \eqref{Delta-u1}, and $\Delta u_{2, \delta} = 0$ in $B_{r_3} \setminus B_{r_*}$, by \eqref{Delta-u2} \footnote{We recall that $r_* = \sqrt{r_2 r_3}$ by \eqref{def-r*}.}, one can represent $u_{\delta}$, $u_{1, \delta}$, and $u_{2, \delta}$ of the forms
\begin{equation}\label{re-u}
u_{\delta} = a_0 + b_0 \ln r + \sum_{n \ge 1} \sum_{\pm} (a_{n, \pm} r^n + b_{n, \pm} r^{-n}) e^{\pm i n \theta} \mbox{ in } B_{r_3} \setminus B_{r_2},
\end{equation}
\begin{equation}\label{re-u1}
u_{1, \delta} = c_0 + d_0 \ln r + \sum_{n \ge 1} \sum_{\pm}(c_{n, \pm} r^n + d_{n, \pm} r^{-n}) e^{\pm i n \theta} \mbox{ in } B_{r_3} \setminus B_{r_2}, 
\end{equation}
and 
\begin{equation}\label{re-u2}
u_{2, \delta} = e_0 + f_0 \ln r + \sum_{n \ge 1} \sum_{\pm}(e_{n, \pm} r^n + f_{n, \pm} r^{-n}) e^{\pm i n \theta} \mbox{ in } B_{r_3} \setminus B_{r_*},
\end{equation}
for $a_0, b_0, c_0, d_0, e_0, f_0, a_{n, \pm}, b_{n, \pm}, c_{n, \pm}, d_{n, \pm}, e_{n, \pm}, f_{n, \pm} \in \mC$ ($n \ge 1$). We derive from \eqref{transmission-r2},  \eqref{re-u}, and \eqref{re-u1} that
\begin{equation*} 
\left\{\begin{array}{ll}
\dsp a_{n, \pm} r_2^n + b_{n, \pm} r_2^{-n}  = & c_{n, \pm} r_2^n + d_{n, \pm} r_2^{-n}, \\[6pt]
\dsp a_{n, \pm} r_2^n - b_{n, \pm} r_2^{-n}  = &(1 - i \delta)  \big( c_{n, \pm} r_2^n - d_{n, \pm} r_2^{-n} \big)
\end{array} \right. \quad \mbox{ for } n \ge 1, 
\end{equation*}
and 
\begin{equation*}
\left\{\begin{array}{rl}
\dsp a_0  + b_0 \ln r_2  = & c_0  + d_0 \ln r_2, \\[6pt]
\dsp b_0  = & (1 - i \delta) d_0. 
\end{array} \right. 
\end{equation*}
This implies
\begin{equation}\label{abcd0} \left\{
\begin{array}{ll}
\dsp a_{n, \pm} = \frac{2 - i \delta}{2} c_{n, \pm} + \frac{i \delta}{2 } d_{n, \pm} r_2^{-2n}, \\[6pt]
\dsp b_{n, \pm} =  \frac{ i \delta}{2} c_{n, \pm} r_2^{2n}+  \frac{2 - i \delta}{2 } d_{n, \pm}, 
\end{array} \right. \quad \mbox{ for } n \ge 1, 
\end{equation}
and 
\begin{equation}\label{abcd0-1}
\left\{\begin{array}{rll}
\dsp a_0    &= & c_0  + i \delta d_0 \ln r_2, \\[6pt]
\dsp b_0  &= & (1 - i \delta) d_0. 
\end{array} \right. 
\end{equation}
Since
\begin{multline*}
u_{\delta} - u_{1, \delta} =  a_0 + b_0 \ln r + \sum_{n \ge 1} \sum_{\pm} (a_{n, \pm} r^n + b_{n, \pm} r^{-n}) e^{\pm i n \theta} \\[6pt] - c_0 - d_0 \ln r - \sum_{n \ge 1} \sum_{\pm}(c_{n, \pm} r^n + d_{n, \pm} r^{-n}) e^{ \pm i n \theta} 
\end{multline*}
in $B_{r_3} \setminus B_{r_2}$, it follows from \eqref{abcd0} and \eqref{abcd0-1} that, in $B_{r_3} \setminus B_{r_2}$, 
\begin{multline}\label{u-u1}
u_{\delta} - u_{1, \delta} =  i \delta d_0 (\ln r_2 - \ln r) -  \frac{i \delta }{2}\sum_{n \ge 1} \sum_{\pm} ( c_{n, \pm} - d_{n, \pm} r_2^{-2n} ) r^n e^{\pm i n \theta} \\[6pt]+  \frac{i \delta }{2}\sum_{n \ge 1} \sum_{\pm} ( c_{n, \pm} r_2^{2n} -   d_{n, \pm} ) r^{-n}e^{ \pm i n \theta}.
\end{multline}
Similarly, we derive from \eqref{transmission-r3},  \eqref{re-u1}, and \eqref{re-u2} that 
\begin{equation*}
\left\{\begin{array}{ll}
\dsp e_{n, \pm} r_3^n + f_{n, \pm} r_3^{-n}  = & c_{n, \pm} r_3^n + d_{n, \pm} r_3^{-n}, \\[6pt]
\dsp e_{n, \pm} r_3^n - f_{n, \pm} r_3^{-n}  = &\dsp  (1 - i \delta)  \big( c_{n, \pm} r_3^n - d_{n, \pm} r_3^{-n} \big),  
\end{array} \right. \quad \mbox{ for } n \ge 1, 
\end{equation*}
and 
\begin{equation*}
\left\{\begin{array}{rll}
\dsp e_0 + f_0 \ln r_3  &= & c_0  + d_0 \ln r_3, \\[6pt]
\dsp f_0  & = &\dsp  (1 - i \delta)  d_0. 
\end{array} \right.
\end{equation*}
This implies
\begin{equation}\label{efcd0} \left\{
\begin{array}{ll}
\dsp e_{n, \pm} = \frac{2 - i \delta}{2} c_{n, \pm} + \frac{i \delta}{2 } d_{n, \pm} r_3^{-2n}, \\[6pt]
\dsp f_{n, \pm} =  \frac{ i \delta}{2} c_{n, \pm} r_3^{2n}+  \frac{2 - i \delta}{2 } d_{n, \pm}, 
\end{array} \right. \quad \mbox{ for } n \ge 1, 
\end{equation}
and 
\begin{equation}\label{efcd0-1}
\left\{\begin{array}{rl}
\dsp e_0    = & c_0  + i \delta d_0 \ln r_3, \\[6pt]
\dsp f_0  = & (1 - i \delta) d_0. 
\end{array} \right. 
\end{equation}
Since
\begin{multline*}
u_{1, \delta} - u_{2, \delta} = c_0 + d_0 \ln r + \sum_{n \ge 1} \sum_{\pm} (c_{n, \pm} r^n + d_{n, \pm} r^{-n}) e^{\pm i n \theta} \\[6pt]
- e_0 - f_0 \ln r - \sum_{n \ge 1} \sum_{\pm}(e_{n, \pm} r^n + f_{n, \pm} r^{-n}) e^{ \pm i n \theta} 
\end{multline*}
in $B_{r_3} \setminus B_{r_*}$, it follows from \eqref{efcd0} and \eqref{efcd0-1} that, in $B_{r_3} \setminus B_{r_*}$, 
\begin{align}\label{u1-u2}
u_{1, \delta} - u_{2, \delta} = & -  i \delta d_0 (\ln r_3 - \ln r) + \frac{i \delta }{2}\sum_{n \ge 1} \sum_{\pm} ( c_{n, \pm} 
- d_{n, \pm} r_3^{-2n} ) r^n e^{ \pm i n \theta} \nonumber\\[6pt]   & -  \frac{i \delta }{2}\sum_{n \ge 1} \sum_{\pm} ( c_{n, \pm} r_3^{2n} -   d_{n, \pm} ) r^{-n}e^{\pm i n \theta}. 
\end{align}
A combination of \eqref{u-u1} and \eqref{u1-u2} yields, in $B_{r_3} \setminus B_{r_*}$, 
\begin{align}\label{u-u2}
u_\delta - u_{2, \delta} = & i \delta d_0 (\ln r_2 - \ln r_3) + \frac{i \delta }{2} \sum_{n \ge 1} \sum_{\pm} c_{n, \pm} (r_2^{2n } - r_3^{2n}) r^{-n} e^{\pm i n \theta} \nonumber \\[6pt]
& +  \frac{i\delta}{2} \sum_{n \ge 1} \sum_{\pm} d_{n, \pm} (r_2^{-2n} - r_3^{-2n}) r^n e^{\pm i n \theta} . 
\end{align}

We now use the removing of localized singularity technique introduced in \cite{Ng-Negative-Cloaking}. Set 
\begin{equation*}
U_\delta = \left\{\begin{array}{cl}
u_\delta  - \hat u_\delta &  \mbox{ if } x \in \Omega \setminus B_{r_*}, \\[6pt] 
u_{2, \delta} & \mbox{ if } x \in B_{r_*}, 
\end{array} \right.
\end{equation*} 
where 
\begin{equation}\label{def-uhat}
\hat u_\delta =  i \delta d_0 (\ln r_2 - \ln r_3)  +  \frac{i \delta}{2} \sum_{n \ge 1} \sum_{\pm} \big( c_n r_2^{2n}- c_n r_3^{2n} \big) r^{-n} e^{\pm i n \theta} \quad  \mbox{ for } |x| \ge r_*. 
\end{equation}
As in \cite{Ng-Negative-Cloaking},  we remove $\hat u_\delta$ from $u_\delta$ in $\Omega \setminus B_{r_*}$. The function $\hat u_\delta$ contains very high modes and creates a trouble for estimating $u_\delta - u_{2, \delta}$ on $\partial B_{r_*}$ (to obtain an estimate for $u_\delta$). However this term can be negligible for large $|x|$ since $r^{-n}$ is small for large r and large n; hence $u_\delta - \hat u_\delta$ well approximates $u_\delta$ for $|x|$ large enough. This is the spirit of the removing of localized singularity technique. 

\medskip

We next estimate 
\begin{equation*}
[U_\delta] \mbox{ and } \big[\hat A \nabla U_\delta \cdot x/|x| \big]  \mbox{ on } \partial B_{r_*}. 
\end{equation*}
Here and in what follows $[U]$ and $\big[\hat A \nabla U_\delta \cdot x/|x| \big]$ denote the jumps of $U_\delta$ and $\hat A \nabla U_\delta \cdot x/|x|$ on $\partial B_{r_*}$. 

\medskip
From \eqref{u-u2} and \eqref{def-uhat}, we have
\begin{equation}\label{U}
[U_\delta] =    \frac{i \delta }{2}\sum_{n \ge 1} \sum_{\pm} d_{n, \pm} ( r_2^{-2n} - r_3^{-2n} )  r_*^n e^{\pm i n \theta} \quad \mbox{ on } \partial B_{r_*}.
\end{equation}
This implies 
\begin{equation}\label{jump-1}
\left\| [U_\delta] \right\|_{H^{1/2}(\partial B_{r_*})}^2 \le C \delta^2 \sum_{n \ge 1} \sum_{\pm} n |d_{n, \pm}|^2  r_2^{-4n} r_*^{2n}. 
\end{equation}
Since, by \eqref{est-F},  
\begin{equation}
\| u_{1, \delta}\|_{H^1(B_{r_3} \setminus B_{r_2})}^2 \le  \frac{C}{ \delta} \|f \|_{L^2} \| u_\delta\|_{L^2(\Omega \setminus B_{r_3})}, 
\end{equation}
and $\Delta u_{1, \delta} = 0 $ in $B_{r_3} \setminus B_{r_2}$, it follows  that 
\begin{equation}\label{u1-r2}
\| u_{1, \delta} \|_{H^{1/2}(\partial B_{r_2})}^2 + \| \partial_r u_{1, \delta} \|_{H^{-1/2}(\partial B_{r_2})}^2 \le  \frac{C}{ \delta} \|f \|_{L^2} \| u_\delta\|_{L^2(\Omega \setminus B_{r_3})}
\end{equation}
and 
\begin{equation}\label{u1-r3}
\| u_{1, \delta} \|_{H^{1/2}(\partial B_{r_3})}^2 + \| \partial_r u_{1, \delta} \|_{H^{-1/2}(\partial B_{r_3})}^2 \le  \frac{C}{ \delta} \|f \|_{L^2} \| u_\delta\|_{L^2(\Omega \setminus B_{r_3})}. 
\end{equation}
A combination of  \eqref{re-u1} and \eqref{u1-r2} yields 
\begin{equation}\label{bdr1}
|c_0|^2 + |d_0|^2 + \sum_{n \ge 1} \sum_{\pm} n\left(|c_{n, \pm}|^2 r_2^{2n} + |d_{n, \pm}|^2 r_2^{-2n} \right) \le \frac{C}{ \delta} \|f \|_{L^2} \| u_\delta\|_{L^2(\Omega \setminus B_{r_3})},  
\end{equation}
and a combination of   \eqref{re-u1} and \eqref{u1-r3} implies
\begin{equation}\label{bdr2}
|c_0|^2 + |d_0|^2 + \sum_{n \ge 1} \sum_{\pm} n \left(|c_{n, \pm}|^2 r_3^{2n} + |d_{n, \pm}|^2 r_3^{-2n} \right) \le \frac{C}{ \delta } \|f \|_{L^2} \| u_\delta\|_{L^2(\Omega \setminus B_{r_3})}. 
\end{equation}
Similarly, 
\begin{equation}\label{bdr3}
|a_0|^2 + |b_0|^2 + \sum_{n \ge 1} \sum_{\pm} n \left(|a_{n, \pm}|^2 r_3^{2n} + |b_{n, \pm}|^2 r_3^{-2n} \right) \le \frac{C}{ \delta } \|f \|_{L^2} \| u_\delta\|_{L^2(\Omega \setminus B_{r_3})}. 
\end{equation}
We derive from \eqref{abcd0}, \eqref{abcd0-1}, \eqref{bdr2}, and \eqref{bdr3} that
\begin{equation}\label{bdr4}
|c_0|^2 + |d_0|^2 + \sum_{n \ge 1} \sum_{\pm} n \left(|c_{n, \pm}|^2 r_3^{2n} + \delta^2 |d_{n, \pm}|^2 r_3^{2n} r_2^{-4n} \right) \le \frac{C}{ \delta } \|f \|_{L^2} \| u_\delta\|_{L^2(\Omega \setminus B_{r_3})}. 
\end{equation}
Since $r_* = \sqrt{r_2 r_3}$, by the H\"older inequality, we have
\begin{equation}\label{ineq-Holder}
\delta^2 \sum_{n \ge 1} \sum_{\pm } n |d_n|^2 r_2^{-4n} r_*^{2n} \le \delta \left( \delta^2\sum_{n \ge 1} \sum_{\pm} n |d_n|^2 r_2^{-4n} r_3^{2n} \right)^{1/2} \left( \sum_{n \ge 1} \sum_{\pm} n |d_n|^2 r_2^{-2n} \right)^{1/2}. 
\end{equation}
A combination of \eqref{jump-1}, \eqref{bdr1}, \eqref{bdr4}, and \eqref{ineq-Holder} yields 
\begin{equation}\label{est-jump1}
\left\| \left[U_\delta\right] \right\|_{H^{1/2}(\partial B_{r_*})}^2 \le  C \|f \|_{L^2} \| u_\delta\|_{L^2(\Omega \setminus B_{r_3})}. 
\end{equation}
Similarly, we have 
\begin{equation}\label{est-jump2}
\left\| \left[\hat A \nabla U_\delta \cdot x/ |x| \right] \right\|_{H^{-1/2}(\partial B_{r_*})}^2 \le  C \|f \|_{L^2} \| u_\delta\|_{L^2(\Omega \setminus B_{r_3})}. 
\end{equation}

On the other hand, from \eqref{def-uhat}, we have 
\begin{equation}\label{est-uhat}
\| \hat u_\delta\|_{H^1(\Omega \setminus B_{r_3})}^2 \le C \big( \delta^2 d_0^2 + \delta^2 \sum_{n \ge 1} \sum_{\pm} n |c_n|^2 r_3^{2n} \big). 
\end{equation}
We derive from  \eqref{bdr1}, \eqref{bdr2}, and \eqref{est-uhat} that 
\begin{equation}\label{est-out}
\| \hat u_\delta\|_{H^1(\Omega \setminus B_{r_3})}^2 \le  C \delta \|f \|_{L^2} \| u_\delta\|_{L^2(\Omega \setminus B_{r_3})},
\end{equation}
which implies, since $U_\delta = u_\delta - \hat u_\delta$ in $\Omega \setminus B_{r_3}$, 
\begin{equation}\label{est-out1}
\| \hat u_\delta\|_{H^1(\Omega \setminus B_{r_3})} \le C \delta^{1/2} \big( \| f\|_{L^2} + \| U_\delta\|_{L^2(\Omega \setminus B_{r_3})} \big). 
\end{equation}
It follows from
 \eqref{est-jump1}, \eqref{est-jump2}, and \eqref{est-out} that
\begin{equation}\label{jump-12}
\left\| \left[U_\delta \right] \right\|_{H^{1/2}(\partial B_{r_*})}^2 + \left\| \left[\hat A \nabla U_\delta \cdot x/ |x| \right] \right\|_{H^{-1/2}(\partial B_{r_*})}^2 \le C  \Big( \|f \|_{L^2} \| U_\delta\|_{L^2(\Omega \setminus B_{r_3})} + \delta^{1/2} \| f\|_{L^2}^2 \Big). 
\end{equation}
Since $\dive(\hat A U_\delta) = f$ in $\Omega \setminus \partial B_{r_*}$, $U_\delta \in H^1(\Omega \setminus \partial B_{r_*})$, and $U_\delta = - \hat u_\delta$ on $\partial \Omega$, we derive from \eqref{est-out1} and \eqref{jump-12} that  
\begin{equation}\label{conclusion}
\| U_\delta \|_{H^1(\Omega \setminus B_{r_*})} + \| U_\delta\|_{H^1(B_{r_*})} \le C \| f\|_{L^2}. 
\end{equation}
A combination of \eqref{est-out1} and \eqref{conclusion} yields 
\begin{equation}\label{est-out2}
\| \hat u_\delta\|_{H^1(\Omega \setminus B_{r_3})} \to 0 \mbox{ as } \delta \to 0. 
\end{equation}
We claim that
\begin{equation}\label{claim1}
\left[U_\delta \right]  \to 0 \mbox{ weakly in } H^{1/2}(\partial B_{r_*}) \; \;  \mbox{ and } \; \;  \left[\hat A \nabla U_\delta \cdot x/ |x| \right]  \to 0 \mbox{ weakly in } H^{-1/2}(\partial B_{r_*}). 
\end{equation}
Assuming \eqref{claim1} holds, we have
\begin{equation*}
U_\delta \to U_0 \mbox{ weakly in } H^1(\Omega \setminus B_{r_3})
\end{equation*}
where $U_0 \in H^1_0(\Omega)$ (by \eqref{est-out2}) is the unique solution to the equation 
\begin{equation*}
\dive(\hat A \nabla U_0) = f \mbox{ in } \Omega. 
\end{equation*}
The conclusion now follows from \eqref{est-out1}. 

\medskip
It remains to prove \eqref{claim1}. We only prove that
\begin{equation*}
\left[U_\delta \right]  \to 0 \mbox{ weakly in } H^{1/2}(\partial B_{r_*}), 
\end{equation*}
the proof of the statement 
\begin{equation*}
\left[\hat A \nabla U_\delta \cdot x/ |x| \right]  \to 0 \mbox{ weakly in } H^{-1/2}(\partial B_{r_*})
\end{equation*}
follows similarly.  Indeed, since $\| U_\delta \|_{H^1(\Omega \setminus B_{r_3})} \le C \| f\|_{L^2}$, it follows from \eqref{est-out1} that 
\begin{equation*}
|a_0|^2 + |b_0|^2 + \sum_{n \ge 1} n \left(|a_{n, \pm}|^2 r_3^{2n} + |b_{n, \pm}|^2 r_3^{-2n} \right) \le C \| f\|_{L^2}.
\end{equation*}
 We derive from \eqref{abcd0} and \eqref{abcd0-1} that 
\begin{equation*}
|d_{n, \pm}| \le C(n) \|f \|_{L^2}, 
\end{equation*}
for some $C(n)$ depending only on $n$, $r_2$, and $r_3$. Since 
\begin{equation*}
\left\| [U_\delta \right] \|_{H^{1/2}(\partial B_{r_*})}  \le C \| f\|_{L^2}, 
\end{equation*}
by \eqref{jump-1}, \eqref{est-out},  and \eqref{conclusion}, it follows from \eqref{U}
that 
\begin{equation*}
\left[U_\delta \right]  \to 0 \mbox{ weakly in } H^{1/2}(\partial B_{r_*}),
\end{equation*}
 
\medskip 
The proof is complete. \proofend

\subsection{Proof of Theorem~\ref{thm1} in the three dimensional  quasistatic regime}\label{part2}

\medskip The proof in the three dimensional quasistatic case follows similarly as the one in two dimensions. We also have  $\Delta u_\delta = \Delta u_{1, \delta} = 0$ in $B_{r_3} \setminus B_{r_2}$, and $\Delta u_{2, \delta} = 0$ in $B_{r_3} \setminus B_{r_*}$. Hence $u_{\delta}$, $u_{1, \delta}$, and $u_{2, \delta}$ can be written under the forms
\begin{equation}\label{re-u-3d}
u_{\delta}(x) = a_0 + \frac{b_0}{r} + \sum_{n=1}^\infty \sum_{l=-n}^n (a_{n, l} |x|^n + b_{n, l} |x|^{-n - 1}) Y^l_n(x/|x|)  \mbox{ in } B_{r_3} \setminus B_{r_2},
\end{equation}
\begin{equation}\label{re-u1-3d}
u_{1, \delta}(x) =  c_0 + \frac{d_0}{r} + \sum_{n=1}^\infty \sum_{l=-n}^n (c_{n, l} |x|^n + d_{n, l} |x|^{-n - 1}) Y^l_n(x/|x|)  \mbox{ in } B_{r_3} \setminus B_{r_2}, 
\end{equation}
and 
\begin{equation}\label{re-u2-3d}
u_{2, \delta}(x) =  e_0 + \frac{f_0}{r} + \sum_{n=1}^\infty \sum_{l=-n}^n (e_{n, l} |x|^n + f_{n, l} |x|^{-n-1}) Y^l_n(x/|x|)  \mbox{ in } B_{r_3} \setminus B_{r_*},
\end{equation}
for some $a_0, b_0, c_0, d_0, e_0, f_0, a_{n, l}, b_{n, l}, c_{n, l}, d_{n, l}, (e_{n, l}, f_{n, l} \in \mC$ ($n \ge 1$, $-n \le l \le n$).
Here $Y_n^l$ is the spherical harmonic function of degree $n$ and of order $l$. The details are left to the reader.  \proofend

\subsection{Proof of Theorem~\ref{thm1} in the finite frequency regime}\label{part3}

The proof in this case is similar to the one in the quasi static case though it is more complicated.  We will present necessary modifications in the two dimensional case. The three dimensional case follows similarly. For notational ease, we will assume $k=1$. 

\medskip
Let $d=2$ and $k=1$.  
Using \eqref{WP} and \eqref{WP1} and applying the same method used in the proof of  \cite[Lemma 1]{Ng-Complementary}, we obtain, for small $\delta$, 
\begin{equation*}
\| u_\delta\|_{H^1(\Omega)}^2 \le C \left(\frac{1}{\delta} \left| \int_\Omega f \bar u_\delta \right| + \|f \|_{L^2}^2 \right).  
\end{equation*}
This implies 
\begin{equation}
\| u_\delta\|_{H^1(\Omega)}^2 \le C \left(\frac{1}{\delta} \| f\|_{L^2} \| u_\delta\|_{L^2(\Omega \setminus B_{r_3})} + \|f \|_{L^2}^2 \right).  
\end{equation}
We  have  
\begin{equation}\label{important}
\Delta u_\delta + k^2 u_\delta = \Delta u_{1, \delta} + k^2 u_{1, \delta} = 0 \mbox{ in } B_{r_3} \setminus B_{r_2}  \mbox{ and } \Delta u_{2, \delta} + k^2 u_{2, \delta}= 0 \mbox{ in } B_{r_3} \setminus B_{r_*}
\end{equation}
by \eqref{F-b}, \eqref{the-choice}, and  \eqref{choice-sigma}.  From \eqref{important}, 
one can represent $u_{\delta}$, $u_{1, \delta}$, and $u_{2, \delta}$ of the forms
\begin{equation}\label{re-u-k}
u_{\delta} =  a_0 \hat J_0(|x|) + b_0 \hat Y_0(|x|) +   \sum_{n=1}^\infty \sum_{\pm} \big[a_{n, \pm} \hat J_{n}(|x|) + b_{n, \pm} \hat Y_{n}(|x|) \big] e^{\pm i n \theta} \quad \mbox{ in } B_{r_3} \setminus B_{r_2},
\end{equation}
\begin{equation}\label{re-u1-k}
u_{1, \delta} =  c_0 \hat J_0(|x|) + d_0 \hat Y_0(|x|)  + \sum_{n=1}^\infty  \big[c_{n, \pm} \hat J_{n}(|x|) + d_{n, \pm} \hat Y_{n}(|x|) \big] e^{\pm i n \theta} \quad \mbox{ in } B_{r_3} \setminus B_{r_2}, 
\end{equation}
and 
\begin{equation}\label{re-u2-k}
u_{2, \delta} =  e_0 \hat J_0(|x|) + f_0 \hat Y_0 (|x|) + \sum_{n=1}^\infty \big[e_{n, \pm} \hat J_{n}(|x|) + f_{n, \pm} \hat Y_{n} (|x|) \big] e^{\pm i n \theta} \quad \mbox{ in } B_{r_3} \setminus B_{r_*},
\end{equation}
for $a_0, b_0, c_0, d_0, e_0, f_0, a_{n, \pm}, b_{n, \pm}, c_{n, \pm}, d_{n, \pm}, e_{n, \pm}, f_{n, \pm} \in \mC$ ($n \ge 1$). Here 
\begin{equation*}
\hat J_n(r) = 2^n n! J_n(r) \quad \mbox{ and } \quad \hat Y_n(r) = -  \frac{\pi}{2^{n} (n-1)!} Y_n(r), 
\end{equation*}
where $J_n$ and $Y_n$ are the Bessel and Neumann functions of order $n$. It follows from \cite[(3.57) and (3.58)]{ColtonKressInverse} that 
\begin{equation}\label{bh1}
\hat J_n (t)  = t^{n}\big[1 + o(1) \big]
\end{equation}
and
\begin{equation}\label{bh2}
\hat Y_n (t)  = t^{-n} \big[1 + o(1) \big], 
\end{equation}
as $n \to + \infty$. Similar to \eqref{abcd0}, we have
\begin{equation}\label{abcd-k}\left\{
\begin{array}{l}
\dsp a_{n, \pm}  = c_{n, \pm} + i \delta c_{n, \pm} \AC_{n, \pm} + i \delta d_{n, \pm} \AD_{n, \pm}, \\[6pt]
\dsp b_{n, \pm} = i \delta c_{n, \pm}  \BC_{n, \pm} + d_{n, \pm} +  i \delta d_{n, \pm} \BD_{n, \pm},
\end{array} \right. \quad \mbox{ for } n \ge 0, 
\end{equation}
and similar to \eqref{efcd0}, we obtain
\begin{equation}\label{efcd-k}\left\{
\begin{array}{l}
\dsp e_{n, \pm}  = c_{n, \pm} + i \delta c_{n, \pm} \EC_{n, \pm} + i \delta d_{n, \pm} \ED_{n, \pm}, \\[6pt]
\dsp f_{n, \pm} = i \delta c_{n, \pm} \FC_{n, \pm} + d_{n, \pm} +  i \delta d_{n, \pm} \FD_{n, \pm}, 
\end{array} \right. \quad \mbox{ for } n \ge 0. 
\end{equation}
Here we denote $a_{0, \pm}, b_{0, \pm}, c_{0, \pm}, d_{0, \pm}, e_{0, \pm}$ and $f_{0, \pm}$ by $a_0/2, b_0/2, c_0/2, d_0/2, e_0/2$, and  $f_0/2$ respectively, and 
\begin{equation}
\AC_{n} = \frac{\hat J_n' \hat Y_n}{\hat J_n \hat Y_n' - \hat J_n' \hat Y_n}(r_2), \quad \AD_{n} =  \frac{\hat Y_n \hat Y_n' }{\hat J_n \hat Y_n' - \hat J_n' \hat Y_n}(r_2), 
\end{equation}
\begin{equation}
\BC_{n} = \frac{\hat J_n \hat J_n'}{\hat Y_n \hat J_n' - \hat Y_n' \hat J_n}(r_2), \quad \BD_n = \frac{\hat Y_n' \hat J_n}{\hat Y_n \hat J_n' - \hat Y_n' \hat J_n}(r_2), 
\end{equation}
\begin{equation*}
\EC_n =  \frac{\hat J_n' \hat Y_n}{\hat J_n \hat Y_n' - \hat J_n' \hat Y_n}(r_3), \quad \ED_n = \frac{ \hat Y_n \hat Y_n'}{\hat J_n \hat Y_n' - \hat J_n' \hat Y_n} (r_3), 
\end{equation*}
and 
\begin{equation*}
\FC_n = \frac{\hat J_n \hat J_n'}{\hat Y_n \hat J_n' - \hat Y_n' \hat J_n} (r_3), \quad \FD_n = \frac{\hat Y_n' \hat J_n}{\hat Y_n \hat J_n' - \hat Y_n' \hat J_n} (r_3). 
\end{equation*}
Then, in $B_{r_3} \setminus B_{r_2}$,  
\begin{multline}\label{u-u1-1}
u_{\delta} - u_{1, \delta} = \sum_{n \ge 0} \sum_{\pm}  i \delta (\AC_n c_{n, \pm} + \AD_n d_{n, \pm}) \hat J_n (|x|) e^{\pm i n \theta}  \\[6pt]
+\sum_{n \ge 0}\sum_{\pm} i \delta (\BC_n c_{n, \pm} + \BD_n d_{n, \pm}) \hat Y_n (|x|) e^{\pm i n \theta} 
\end{multline}
and, in $B_{r_3} \setminus B_{r_*}$, 
\begin{multline}\label{u1-u2-1}
u_{1, \delta} - u_{2, \delta} =   - \sum_{n \ge 0} \sum_{\pm} i \delta (\EC_n c_{n, \pm} + \ED_n d_{n, \pm}) \hat J_n (|x|) e^{\pm i n \theta} \\[6pt]
 - \sum_{n \ge 0} \sum_{\pm} i \delta (\FC_n c_{n, \pm} + \FD_n d_{n, \pm}) \hat Y_n (|x|) e^{\pm i n \theta}.
\end{multline}
A combination of \eqref{u-u1-1} and \eqref{u1-u2-1} yields, in $B_{r_3} \setminus B_{r_*}$,  
\begin{align}\label{u-u2-1}
u_\delta - u_{2, \delta} = & \sum_{n \ge 0} \sum_{\pm} i \delta \Big[ c_{n, \pm} (\AC_n - EC_n)  +  d_{n, \pm} (\AD_n - ED_n) \Big] \hat J_n (|x|) e^{\pm i n \theta} \nonumber\\[6pt]
& + \sum_{n \ge 0} \sum_{\pm} i \delta   \Big[ (\BC_n - \FC_n) c_{n, \pm} + ( \BD_n - \FD_n) d_{n, \pm} \Big]  \hat Y_n(|x|) e^{\pm i n \theta}.
\end{align}

We now use the removing of localized singularity technique as in the quasistatic case. Set
\begin{equation*}
\hat u_\delta (x)=  \sum_{n \ge 0} \sum_{\pm} i \delta   \Big[ (\BC_n - \FC_n) c_{n, \pm} + ( \BD_n - \FD_n) d_{n, \pm} \Big]  \hat Y_n(|x|) e^{\pm i n \theta}.
\end{equation*}
and define
\begin{equation*}
U_\delta = \left\{\begin{array}{cl} u_\delta - \hat u_\delta & \mbox{ if } x \in \Omega \setminus B_{r_*}, \\[6pt]
u_{2, \delta} &  \mbox{ if } x \in B_{r_*}. 
\end{array} \right.
\end{equation*}
Using \eqref{bh1} and \eqref{bh2}, we have
\begin{equation}\label{ACD}
\AC_n  = -  \frac{1}{2}\big[1 + o(1)\big], \quad \AD_n = \frac{1}{2} r_2^{-2n} \big[1 + o(1)\big], 
\end{equation}
and 
\begin{equation}\label{BCD}
 \BC_n = \frac{1}{2} r_2^{2n} \big[1 + o(1)\big], \quad \BD_n = - \frac{1}{2} \big[1 + o(1)\big].
\end{equation}
Similarly, we obtain
\begin{equation}\label{ECD}
\EC_n  =  - \frac{1}{2}\big[1 + o(1)\big], \quad \ED_n = \frac{1}{2} r_3^{-2n} \big[1 + o(1)\big], 
\end{equation}
and 
\begin{equation}\label{FCD}
\FC_n = \frac{1}{2} r_3^{2n} \big[1 + o(1)\big], \quad \FD_n = - \frac{1}{2}\big[1 + o(1)\big].
\end{equation}
Since (see, e.g., \cite[(3.56)]{ColtonKressInverse}) 
\begin{equation*}
\hat J_n' (r) \hat Y_n(r) - \hat J_n (r) \hat Y_n' (r) = C_n r^{-1}, 
\end{equation*}
it follows that 
\begin{equation}\label{bh3}
|c_{n, \pm}|^2 + |d_{n, \pm}|^2 \le C_{n, r} \left( |c_{n, \pm} \hat J_n(r) + d_{n, \pm} \hat Y_n(r)|^2 + |c_{n, \pm} \hat J_n'(r) + d_{n, \pm} \hat Y_n'(r)|^2 \right). 
\end{equation}
Combining  \eqref{bh1}, \eqref{bh2}, and \eqref{bh3}, as in \eqref{bdr2}, we obtain 
 \begin{equation}\label{bdr2-k}
\sum_{n \ge 0} \sum_{\pm}(n+1) \left(|c_{n, \pm}|^2 r_3^{2n} + |d_{n, \pm}|^2 r_3^{-2n} \right) \le  C \left(\frac{1}{\delta} \| u_\delta\|_{L^2(\Omega \setminus B_{r_3})} \| f\|_{L^2} + \|f \|_{L^2}^2 \right)
\end{equation}
and 
\begin{equation}\label{bdr2-k-1}
\sum_{n \ge 0} \sum_{\pm} (n+1) \left(|c_{n, \pm}|^2 r_2^{2n} + |d_{n, \pm}|^2 r_2^{-2n} \right) \le  C \left(\frac{1}{\delta} \| u_\delta\|_{L^2(\Omega \setminus B_{r_3})} \| f\|_{L^2} + \|f \|_{L^2}^2 \right)
\end{equation}
Similarly, 
\begin{equation}\label{bdr3-k}
\sum_{n \ge 0} \sum_{\pm} (n+1) \left(|a_{n, \pm}|^2 r_3^{2n} + |b_{n, \pm}|^2 r_3^{-2n} \right) \le  C \left(\frac{1}{\delta} \| u_\delta\|_{L^2(\Omega \setminus B_{r_3})} \| f\|_{L^2} + \|f \|_{L^2}^2 \right).   
\end{equation}
We derive from \eqref{bdr2-k} and \eqref{bdr3-k} that 
\begin{equation}\label{bdr4-k}
\sum_{n \ge n_0} \sum_{\pm} (n+1) \left(|c_{n, \pm}|^2 r_3^{2n} + \delta^2 |d_{n, \pm}|^2 r_3^{2n} r_2^{-4n} \right) \le  C \left(\frac{1}{\delta} \| u_\delta\|_{L^2(\Omega \setminus B_{r_3})} \| f\|_{L^2} + \|f \|_{L^2}^2 \right), 
\end{equation}
for some $n_0$ large enough. 

\medskip
As in the proof of Theorem~\ref{thm1} in the quasistatic case,  using \eqref{bdr2-k-1} and \eqref{bdr4-k}, we have 
\begin{equation*}
\big\| [U_\delta] \big\|_{H^{1/2}(\partial B_{r_*})}^2 + \big\| [\hat A \nabla U_\delta \cdot \eta] \big\|_{H^{-1/2}(\partial B_{r_*})}^2  \le C \big( \| f\|_{L^2} \| u_\delta\|_{L^2(\Omega \setminus B_{r_3})} + \delta \| f\|_{L^2}^2 \big). 
\end{equation*}

\medskip
The proof is now similar to the one in the quasistatic case. The uniqueness of the limit of $U_\delta$ follows from \eqref{WP}. The details are left to the reader. \proofend

\section{Other constructions of superlenses}\label{sect-comments}

The construction of the superlens given by \eqref{defA} and \eqref{def-ss} is not restricted to the Kelvin transform $F$ w.r.t. $\partial B_{r_2}$. In fact, 
using the study of reflecting complementary media in \cite{Ng-Complementary}, we can extend this construction further. We confine ourselves to a class of radial reflections for which the formulae for $A$ and $\Sigma$ are explicit even though  general reflections as in \cite{Ng-Complementary} can be used. 

\medskip
Fix $\alpha, \beta > 1$ such that \footnote{One can choose $\alpha$ and $\beta$ such that $\alpha \beta - \alpha - \beta \ge 0$. However, these expression of $A_1$ and $\hat A_1$ below are more involved in this case.}
\begin{equation}\label{ab}
\alpha \beta - \alpha - \beta = 0.
\end{equation}
Let $F_1: B_{r_2} \setminus \{ 0 \} \to \mR^d \setminus \bar B_{r_2}$ and $G_1: \mR^d \setminus \bar B_{r_3} \to B_{r_3} \setminus \{0\}$ be defined as follows: 
\begin{equation*}
 F_1(x) = r_2^\alpha x/ |x|^\alpha \quad \mbox{ and }  \quad G_1(x) = r_3^\beta x/ |x|^\beta. 
\end{equation*}
Here, $r_1$, $r_2$, and $r_3$ are chosen such that 
\begin{equation*}
r_3/r_1 = r_2^{\alpha}/ r_1^\alpha = m \quad \mbox{ and } \quad \sqrt{r_2 r_3} = m r_0; 
\end{equation*}
which yields
\begin{equation}\label{choice-123}
r_1 = r_0 m^{\frac{\alpha-1}{2\alpha}}, \quad   r_2 = r_0  m^{\frac{\alpha+1}{2\alpha}}, \quad \mbox{ and } \quad r_3 = r_0 m^{\frac{3\alpha-1}{2 \alpha}}. 
\end{equation}

It follows from \eqref{ab} and \eqref{choice-123} that  $G_1 \circ F_1 : B_{r_1} \to B_{r_3}$ satisfies
\begin{equation*}
G_1 \circ F_1 (x) =    m x.
\end{equation*}
Define
\begin{equation}\label{defAtilde}
A_1, \Sigma_1 = \left\{ \begin{array}{cl} 
\big(F_1^{-1}\big)_*I, \; \big(F_1^{-1} \big)_*1 & \mbox{ in } B_{r_2} \setminus B_{r_1}, \\[6pt]
\big(F_1^{-1} \big)_*\big(G_1^{-1} \big)_* I, \;  \big(F_1^{-1} \big)_*\big( G_1^{-1} \big)_* I & \mbox{ in } B_{r_1} \setminus B_{r_0},\\[6pt]
a, \; \sigma & \mbox{ in } B_{r_0}, \\[6pt]
I, \; 1 & \mbox{ otherwise}, 
\end{array} \right. 
\end{equation} 
and 
\begin{equation*} \label{defAtildehat}
\hat A_1, \hat \Sigma_1 = \left\{ \begin{array}{cl} I, 1 & \mbox{ in } \Omega \setminus B_{m r_0}, \\[6pt] 
(G_1)* (F_1)_* a, \;  (G_1)_*  (F_1)_* 1  =  m^{2-d} a(x/ m), \; m^{-d} \sigma(x/m)& \mbox{ in } B_{m r_0}.  
\end{array} \right.
\end{equation*}
One can verify that, in $B_{r_2} \setminus B_{r_1}$,  
\begin{equation}\label{id2}
A_1, \; \Sigma_1  = \frac{r_2^\alpha}{ r^{\alpha}}  \left[ \frac{1}{\alpha - 1} e_{r} \otimes e_{r } + (\alpha -1) \Big( e_{\theta} \otimes e_{\theta}  + e_{\theta} \otimes e_{\varphi}  \Big) \right],  \; 
 (\alpha -1)\frac{r_2^{3 \alpha}}{r^{3 \alpha}}  \mbox{ if } d = 3, 
\end{equation}
and 
\begin{equation}\label{id3}
A_1, \;  \Sigma_1  = \frac{1}{\alpha -1} e_{r} \otimes e_{r } + (\alpha -1) e_{\theta} \otimes e_{\theta} ,  \; (\alpha - 1)\frac{r_2^{2 \alpha}}{r^{2 \alpha}}  \mbox{ if } d =2. 
\end{equation}
and, in $B_{r_1} \setminus B_{r_0}$, 
\begin{equation}\label{id4}
 A_1,  \; \Sigma_1  = m^{d-2} I, \;  m^{d}. 
\end{equation}

We will assume that \eqref{WP}  holds for $(\hat A_1, \hat \Sigma_1)$ instead of $(\hat A,  \hat \Sigma)$
and  
\begin{equation}\label{WP1-1}
\mbox{ equation } \Delta v + k^2 v = 0 \mbox{ in } \Omega \setminus B_{r_2} \mbox{ has only zero solution in } H^1_0(\Omega \setminus B_{r_2}). 
\end{equation}

The following result is a generalization of Theorem~\ref{thm1}. 

\begin{theorem} \label{thm2}
Let  $d=2, \, 3$, $f \in L^{2}(\Omega)$ with $\supp f \subset \Omega \setminus B_{r_{3}}$ and let  $u, \; u_\delta \in  H^1_{0}(\Omega)$ be respectively the unique solutions to 
\begin{equation*}
\dive (s_\delta  A_1 \nabla u_\delta) + s_0 k^2 \Sigma_1 u_{\delta} = f \mbox{ in } \Omega
\end{equation*}
and 
\begin{equation*}
\dive (\hat A_1 \nabla u) + k^2 \hat \Sigma_1 u = f \mbox{ in } \Omega
\end{equation*}
We have 
\begin{equation}\label{key-point-*}
u_{\delta} \to u \mbox{ weakly in } H^{1} (\Omega \setminus B_{r_3}) \mbox{ as } \delta \to 0. 
\end{equation}
\end{theorem}

By taking $\alpha = \beta =2$, we obtain Theorem~\ref{thm1} from Theorem~\ref{thm2}. 

\begin{remark}\label{rem-thickness}
We have $\beta = \alpha/ (\alpha-1)$ by \eqref{ab}. Letting $\alpha \to 1_+$, we derive from \eqref{choice-123} that 
\begin{equation}\label{thickness}
r_1 \to r_0 \quad  \mbox{ and } \quad r_3 \to m r_0. 
\end{equation}
Thus for any $\eps > 0$, there exists a construction such that the magnification of $m$ times for an  object in $B_{r_0}$  takes place for any $\supp f \subset \Omega \setminus B_{mr_0 + \eps}$. 
\end{remark}


\noindent {\bf Proof.} We have
\begin{equation*}
(F_1)_* A_1 = I \mbox{ in } B_{r_3} \setminus B_{r_2} \quad \mbox{ and } \quad  (G_1)_* (F_1)_* A_1 = I \mbox{ in } B_{r_3} \setminus B_{r_*},  
\end{equation*}
and
\begin{equation*}
(F_1)_* \Sigma_1 = 1 \mbox{ in } B_{r_3} \setminus B_{r_2} \quad \mbox{ and } \quad  (G_1)_* (F_1)_* \Sigma_1 = 1 \mbox{ in } B_{r_3} \setminus B_{r_*},  
\end{equation*}
by the definition of $(A_1, \Sigma_1)$. We recall that $r_* = \sqrt{r_2 r_3} = m r_0$ by \eqref{choice-123}.  This implies, by \cite[Lemma 4]{Ng-Complementary},
\begin{equation*}
\Delta u_{1, \delta} + k^2 u_{1, \delta}  =  0 \mbox{ in } B_{r_3} \setminus B_{r_2}  
\end{equation*} 
and
\begin{equation*}
\Delta u_{2, \delta} + k^2 u_{2, \delta} = 0 \mbox{ in } B_{r_3} \setminus B_{r_*}. 
\end{equation*} 
Here, as in the proof of Theorem~\ref{thm1}, we define
\begin{equation*}
u_{1, \delta } = u \circ  F_1^{-1} \mbox{ in } \mR^d \setminus B_{r_3} \quad \mbox{ and } \quad u_{2, \delta} = u_{1, \delta} \circ  G_1^{-1} \mbox{ in } B_{r_3}. 
\end{equation*}
Similar to \eqref{transmission-r2} and \eqref{transmission-r3}, by  \cite[Lemma 4]{Ng-Complementary}, we obtain 
\begin{equation*}
u_{1, \delta} = u_{\delta} \Big|_{+}   \mbox{ on } \partial B_{r_{2}} \quad \mbox{ and } \quad  (1 - i \delta ) A_1 \nabla u_{1, \delta} \cdot \eta = A_1 \nabla u_{ \delta} \cdot \eta  \Big|_{+} \mbox{ on } \partial B_{r_{2}}. 
\end{equation*}
and 
\begin{equation*}
u_{1, \delta} = u_{2, \delta} \mbox{ on } \partial B_{r_{3}} \quad  \mbox{ and }  \quad (1 - i \delta )\partial_{\eta} u_{1, \delta} \Big|_{-}= \partial_{\eta} u_{2, \delta} \mbox{ on } \partial B_{r_{3}}.
\end{equation*}
This proof  is now similar to the one of Theorem~\ref{thm1}. The details are left to the reader. \proofend

\bigskip

\providecommand{\bysame}{\leavevmode\hbox to3em{\hrulefill}\thinspace}
\providecommand{\MR}{\relax\ifhmode\unskip\space\fi MR }
\providecommand{\MRhref}[2]{%
  \href{http://www.ams.org/mathscinet-getitem?mr=#1}{#2}
}
\providecommand{\href}[2]{#2}


\end{document}